\documentclass{amsart}

\usepackage{amssymb,amsmath}

\numberwithin{equation}{section}

\newtheorem{theorem}{Theorem}
\newtheorem{lemma}{Lemma}[section]

\newtheorem{prop}{Proposition}
\newtheorem{remark}{Remark}[section]

\newtheorem{define}{Definition}[section]

\newcommand\pd{{\partial}}

\newcommand{\R}{{\mathbb{R}}}
\newcommand{\N}{{\mathbb{N}}}
\newcommand{\cZ}{{\mathbb{Z}}}
\newcommand{\h}{\mathcal{H}}
\newcommand{\cl}{\mathcal{L}}

\newcommand{\p}{\mathcal{P}}

\newcommand\cT{{\mathcal T}}
\newcommand\T{{\mathbb T}}

\newcommand\ep{\epsilon}

\newcommand\zinf{z\rightarrow+\infty}

\newcommand\limZ{\lim\limits_{Z\rightarrow+\infty}}
\newcommand\limz{\lim\limits_{z\rightarrow+\infty}}

\begin{document}

\title[]
{On the ill-posedness of the Prandtl equations in three space dimensions}

\author{Cheng-Jie Liu}
\address{Cheng-Jie Liu
\newline\indent
Department of Mathematics, City University of Hong Kong,
Hong Kong, P. R. China}
\email{cjliusjtu@gmail.com}

\author{Ya-Guang Wang}
\address{Ya-Guang Wang
\newline\indent
Department of Mathematics, and MOE-LSC, Shanghai Jiao Tong University
\newline\indent
Shanghai, 200240, P. R. China}
\email{ygwang@sjtu.edu.cn}

\author{Tong Yang}
\address{Tong Yang
\newline\indent
Department of Mathematics, Shanghai Jiao Tong University
\newline\indent
Shanghai, 200240, P. R. China
\newline\indent And,
Department of mathematics, City University of Hong Kong
\newline\indent
Hong Kong, P. R. China}
\email{matyang@cityu.edu.hk}


\subjclass[2000]{35M13, 35Q35, 76D10, 76D03, 76N20}

\date{}

\keywords{three-dimensional Prandtl equations, ill-posedness, nonlinear instability, shear flow,
monotonic velocity fields.}

\begin{abstract}
In this paper, we give an instability criterion for the Prandtl equations in three space variables, which shows that the monotonicity condition of tangential velocity fields is not sufficient
for the well-posedness of the three dimensional Prandtl equations, in contrast to the classical well-posedness theory of the Prandtl equations in two space variables under the Oleinik monotonicity assumption of the tangential velocity. Both of linear stability and nonlinear stability are considered. This criterion shows that the monotonic shear flow is linear stable for the three dimensional Prandtl equations if and only if the tangential velocity field direction is invariant with respect to the normal variable, and this result is an exact complement to our recent work \cite{LWY} on the well-posedness theory for the three dimensional Prandtl equations with  special structure.

\end{abstract}

\maketitle

\tableofcontents


\section{Introduction}

The inviscid limit of the viscous flow has been known as a challenging
mathematical problem that contains many unsolved problems. For the incompressible Navier-Stokes equations confined in a domain with  boundary, in particular with the non-slip boundary condition, the justification of the inviscid limit remains basically open, c.f. \cite{e-1} and references therein. The main obstruction
comes from the formation of boundary layers near the physical boundary, in which the tangential velocity component changes dramatically.

The foundation of the boundary layer theories was established by Prandtl \cite{Pra}
in 1904 when he introduced the classical Prandtl equations by considering the
incompressible Navier-Stokes equations with non-slip boundary condition. His
observation reveals that outside the layer of thickness of $\sqrt{\nu}$ with
$\nu$ being the viscosity coefficient, the convection dominates so that the flow
can be described approximately by the incompressible Euler equations, however, within the
layer of thickness of $\sqrt{\nu}$ in the vicinity of the boundary, the
convection and viscosity balance so that the flow is governed by the Prandtl
equations that is degenerate and mixed type.
Since then, there have been a lot of mathematical studies on the Prandtl equations,
however, the existing theories are basically limited to the two space dimensional
case except the one in analytic framework by Sammartino and Caflisch \cite{S-C} and others \cite{Z-Z}. On the other hand, in
two dimensional space, the classical work by Oleinik and her collaborators \cite{Ole}
gives the local in time well-posedness when the tangential velocity
component is
monotone in the normal direction,
by using the  Crocco transformation.
Recently, this well-posedness result of the two dimensional Prandtl equations is re-studied in \cite{AWXY,M-W} by direct energy method. In addition to the monotone condition on the
velocity, if a favorable pressure condition is imposed, then global in time
weak solution was also obtained in two dimensional space, see \cite{X-Z}.

The stability mechanism of the three dimensional Prandtl equations is very challenging and delicate mainly due to the possible appearance of secondary flows in the three dimensional boundary layer flow as explained  in Moore \cite{moore}, and it is an open question proposed by Oleinik and Samokhin in the monograph \cite{Ole}.
Recently, in \cite{LWY} the authors  construct a local solution to the three dimensional Prandtl equations when the tangential velocity field direction is invariant with respect to the normal variable under certain monotonicity condition.
In addition,   this special boundary layer flow is linearly stable with respect to any perturbation, and the global in time weak solution
is also obtained under an additional favorable pressure condition \cite{LWY1}.

The purpose of this paper is to investigate the instability of boundary layer flows in three
space dimensions without the special structure proposed in \cite{LWY}, even when the two tangential velocity components of the background state are monotonic.
 This reveals the essential
difference of the Prandtl equations
between two and three space dimensions. For this, let us first review the recent extensive studies on the
instability of the two space dimensional flow around a background state of shear flow
with non-monotonicity.

In fact, without the monotonicity assumption on the tangential component of
the velocity, boundary separation will occur. For this, there are many
physical observations and mathematical studies. For example,
Van Dommelen and Shen in \cite{van} illustrated the ``Van  Dommelen singularity'' by
considering an impulsively started circular cylinder to show
the blowup of the normal velocity, and
 E and Enquist in \cite{e-2} precisely constructed some finite time blowup
solutions to the two-dimensional Prandtl equations.
Started by Grenier's work in 2000,
there are some extensive investigation on the
instability of the two-dimensional Prandtl equations when the background shear flow
has some degeneracy. Precisely, corresponding to the well known
Rayleigh criterion for the Euler flow,  Grenier \cite{grenier} showed that
 the unstable Euler shear flow  yields instability of the Prandtl equations.
It was shown in \cite{GV-D}
 that a non-degenerate
critical point in  the shear flow of the Prandtl equations leads to
  a strong linear ill-posedness of the Prandtl equations in the
Sobolev space framework. Moreover, \cite{GV-N} strengthens
the result of \cite{GV-D}  for any unstable
shear flow. Along this direction, the ill-posedness in the nonlinear
setting was proved in \cite{guo} to show that
the  Prandtl equations are ill-posed near non-stationary and non-monotonic shear flows so that
the asymptotic boundary layer
expansion is not valid for non-monotonic shear layer flows in Sobolev spaces.

To describe the problem to be studied in this paper,
consider the following incompressible Navier-Stokes equations
\begin{equation}\label{ns}
\begin{cases}
\partial_t {\bf u}^\nu+({\bf u}^\nu\cdot
\nabla){\bf u}^\nu+\nabla p^\nu-\nu\Delta {\bf u}^\nu=0, \\
\nabla\cdot {\bf u}^\nu=0, \\
{\bf u}^\nu|_{z=0}=0,
\end{cases}
\end{equation}
 in $\{t>0, (x, y)\in \T^2, z\in \R^+\}$ with boundary at $\{ z=0\}$, here
${\bf u^\nu}=(u^\nu,v^\nu,w^\nu)^T$. According to
Prandtl's observation, set the ansatz for
${\bf u}^\nu$ near $\{z=0\}$ as
\begin{equation}
\begin{cases}
u^\nu(t,x,y,z)=u(t,x,y,\frac{z}{\sqrt{\nu}})+o(1),\\
v^\nu(t,x,y,z)=v(t,x,y,\frac{z}{\sqrt{\nu}})+o(1),\\
 w^\nu(t,x,y,z)=\sqrt{\nu}w(t,x,y,\frac{z}{\sqrt{\nu}})+o(\sqrt{\nu}),
\end{cases}
\end{equation}
and plug it in the Navier-Stokes equations \eqref{ns}, one finds that the boundary layer profile $(u, v, w)(t,x,y,z)$ (here we replace $\frac{z}{\sqrt{\nu}}$ by $z$ for simplicity of notations) satisfies
\begin{equation}\label{3dpr}
\begin{cases}
\partial_t u+(u\partial_x+v\partial_y+w\partial_z)u+\partial_xp^E(t,x,y,0)=\partial_z^2u,\\
\partial_t v+(u\partial_x+v\partial_y+w\partial_z)v+\partial_yp^E(t,x,y,0)=\partial_z^2v,\\
\partial_xu+\partial_yv+\partial_z w=0,\\
(u, v, w)|_{z=0}=0, \qquad \lim\limits_{z\to+\infty}(u, v)=(u^E, v^E)(t,x,y,0),
\end{cases}
\end{equation}
which is the famous Prandtl layer equations.
Here, the pressure $p^E$ is related to the outer Euler flow
 ${\bf u}^E=(u^E,v^E,0)(t,x,y,0)$ through
$$
\partial_t {\bf u}^E+({\bf u}^E\cdot \nabla){\bf u}^E +\nabla p^E=0.
$$

The main results of this paper show that when the background state is
a shear flow $(u^s(t,z),v^s(t,z),0)$ of \eqref{3dpr} with initial data $(U_s(z),V_s(z))$, even under
the  monotonicity condition that
$U'_s(z), \, V'_s(z)>0$, the Prandtl equations \eqref{3dpr} are both linearly and
nonlinearly unstable under a very general assumption that
\begin{equation}\label{con-1}
\exists ~z_0>0,~s.t.~\frac{d}{dz}(\frac{V'_s}{U'_s})(z_0)\neq 0\quad{\rm or} \quad\frac{d}{dz}(\frac{U'_s}{V'_s})(z_0)\neq 0.
\end{equation}
Note that in \cite{LWY}, existence of solutions to the three space
dimensional Prandtl equations was proved with special structure and in that case,
the tangential components of the solution $(u,v,w)(t,x,y,z)$
satisfies $\frac{\partial}{\partial z}(\frac{v}{u})\equiv 0$. In fact, in this case, the
appearance of the secondary flow that is the key factor in instability is avoided. Thus, by combining with the results obtained in \cite{LWY}, we know that the condition \eqref{con-1} is not only sufficient but also necessary for the linear instability of the three dimensional Prandtl equation \eqref{3dpr} linearized
around the monotonic shear flow $(u^s(t,z), v^s(t,z),0)$.

The rest of the paper will be organized as follows. In Section 2, we first state the main
results on the linear and nonlinear instability of the three-dimensional Prandtl equations with
background state as monotonic shear flow.
Then, we prove the linear instability result of the shear flow in Section 3, and the
nonlinear instability will be studied in Section 4. In the Appendix, we give a well-posedness result for the linearized three-dimensional Prandtl equations in the analytic setting with respect to only one horizontal variable, under the assumption that one component of the tangential velocity field of
background shear flow is monotonic.


\section{Main results}

By assuming that the outer Euler flow is uniform in \eqref{3dpr}, consider the following boundary value problem of three dimensional Prandtl equations in $\Omega\triangleq\{(t,x,y,z):t>0,(x,y)\in\T^2,z\in\R^+\},$
\begin{equation}\label{3dpd}
\begin{cases}
\pd_t u+(u\pd_x+v\pd_y+w\pd_z)u-\pd_z^2u=0,\\
\pd_t v+(u\pd_x+v\pd_y+w\pd_z)v-\pd_z^2v=0,\\
\pd_xu+\pd_yv+\pd_zw=0,\\
(u,v,w)|_{z=0}=0,\qquad \limz (u,v)=(U_0,V_0)
\end{cases}\end{equation}
for positive constants $U_0$ and $V_0$.
 To understand this problem, we start with the simple situation of shear flow.
Let $u^s(t,z)$ and $v^s(t,z)$ be smooth solutions of the heat equations:
\begin{equation}\label{shear}\begin{cases}
\pd_t u^s-\pd_z^2u^s=0,\qquad
\pd_t v^s-\pd_z^2v^s=0,\\
(u^s,v^s)|_{z=0}=0,\qquad \limz (u^s,v^s)=(U_0,V_0),\\
(u^s,v^s)|_{t=0}=(U_s,V_s)(z),
\end{cases}\end{equation}
with $(u^s-U_0,v^s-V_0)$ rapidly tending to $0$ when $\zinf$. It is
straightforward to check that the shear velocity profile $(u^s,v^s,0)(t,z)$ satisfies the problem \eqref{3dpd}.

The question we answer in this paper is
whether such trivial profile is stable even when $u^s(t,z)$ and $v^s(t,z)$ are strictly monotonic in $z>0$.
 For this, we first focus on the linear stability problem, and
 consider the linearization of the problem \eqref{3dpd} around $(u^s,v^s,0)$:
\begin{equation}\label{lin}\begin{cases}
\pd_t u+(u^s\pd_x+v^s\pd_y)u+wu_z^s-\pd_z^2u=0,\quad{\rm in}\quad \Omega,\\
\pd_t v+(u^s\pd_x+v^s\pd_y)v+wv_z^s-\pd_z^2v=0,\quad{\rm in}\quad \Omega,\\
\pd_xu+\pd_yv+\pd_zw=0,\quad{\rm in}\quad \Omega,\\
(u,v,w)|_{z=0}=0,\qquad\limz (u,v)=0.
\end{cases}\end{equation}

To present the linear instability result that is motivated by the work \cite{GV-D} for two dimensional problem, we first introduce
some notations. As in \cite{GV-D}, for any $\alpha, \, m\geq0,$ denote by
\[\begin{split}
L_\alpha^2(\R^+)~&:=~\{f=f(z),z\in\R^+;~\|f\|_{L_\alpha^2}\triangleq\|e^{\alpha z}f\|_{L^2}<\infty\},\\
H_\alpha^m(\R^+)~&:=~\{f=f(z),z\in\R^+;~\|f\|_{H_\alpha^m}\triangleq\|e^{\alpha z}f\|_{H^m}<\infty\}, \\
W_\alpha^{m,\infty}(\R^+)~&:=~\{f=f(z),z\in\R^+;~\|f\|_{W_\alpha^{m,\infty}}\triangleq\|e^{\alpha z}f\|_{W^{m,\infty}}<\infty\},
\end{split}\]
and the functional space for $\forall\beta>0$:
\begin{equation*}\begin{split}
E_{\alpha,\beta}~&:=~\Big\{f=f(x,y,z)=\sum\limits_{k_1,k_2\in\cZ}e^{i(k_1x+k_2y)}f_{k_1,k_2}(z),~\|f_{k_1,k_2}\|_{L_\alpha^2}\leq C_{\alpha,\beta}e^{-\beta\sqrt{k_1^2+k_2^2}}\Big\},\\
\end{split}\end{equation*}
with
\[\|f\|_{E_{\alpha,\beta}}\triangleq\sup\limits_{k_1,k_2\in\cZ}e^{\beta\sqrt{k_1^2+k_2^2}}\|f_{k_1,k_2}\|_{L_\alpha^2}.\]
The same notations are also used for the vector functions without confusion.

As in \cite{GV-D}, we first have the following existence
 result for the problem \eqref{lin} when the data are analytic in the tangential variables $(x,y)$.

\begin{prop}\label{prop_1}
Let $(u^s-U_0,v^s-V_0)\in C\big(\R^+;W_\alpha^{1,\infty}(\R^+)\big).$ Then, there exists a $\rho>0$ such that for all $T$ with $\beta-\rho T>0$, and $(u_0,v_0)\in E_{\alpha,\beta}$, the linear problem \eqref{lin} with the initial data $(u,v)|_{t=0}=(u_0,v_0)$ has a unique solution
\[
(u,v)\in C\Big([0,T);E_{\alpha,\beta-\rho T}\Big),\quad (u,v)(t,\cdot)\in E_{\alpha,\beta-\rho t}.
\]
\end{prop}
The proof of this Proposition is the same as that given in \cite[Proposition 1]{GV-D}, so we omit it for brevity.

If we impose monotonic condition on the tangential velocity components of the shear flow
$(u^s,v^s,0)$, then another well-posedness result can be obtained. For this, similar to the previous notations, we
 introduce the following function spaces: for any $ \alpha,\, \beta>0$, set
\begin{equation}\label{def_H}\begin{split}
K^m_{\alpha}~&:=~\Big\{f=f(x,z)
;~\|f\|_{K^m_\alpha}\triangleq \|e^{\alpha z}f\|_{ H^m(\T_x;L^{2}(\R^+_z))}<\infty\Big\},
\end{split}\end{equation}
and
\begin{equation*}\begin{split}
F^m_{\alpha,\beta}~&:=~\Big\{f=f(x,y,z)=\sum\limits_{k\in\cZ}e^{iky}f_{k}(x,z);\quad
\|f_k\|_{K_\alpha^m}\leq C_{\alpha,\beta}e^{-\beta|k|},~\forall k\Big\}
\end{split}\end{equation*}
with
\[\|f\|_{F^m_{\alpha,\beta}}~\triangleq~\sup\limits_{k\in\cZ}e^{\beta|k|}
\|f_k\|_{K_\alpha^m}.\]
The following result shows that the linear problem \eqref{lin} is still well-posed when the analyticity with respect to one horizontal variable given in Proposition \ref{prop_1} is replaced by some monotonicity  assumption.

\begin{prop}\label{prop_2}
Let $(u^s-U_0,v^s-V_0)\in C(\R^+;W_\alpha^{3,\infty}(\R^+)),\alpha>0$ 
satisfying $u_z^s>0$ and
\[\frac{u_{zz}^s}{u_z^s},~\frac{v_z^s}{u_z^s}\in C\big(\R^+;W^{1,\infty}(\R^+)\big),\]
and the initial data $(u,v)|_{t=0}=(u_0,v_0)$ of the problem \eqref{lin} satisfying
 \begin{equation}\label{ass_ini}
(u_0,v_0)(x,y,z)\in F_{\alpha,\beta}^{m}, \qquad  \pd_z\big(\frac{u_0}{u_z^s(0,\cdot)}\big)~\in~F^{m}_{\alpha,\beta}.
 \end{equation}
 Then, there exists a $\rho>0$ such that for all $T$ with $\beta-\rho T>0$,  the linear problem \eqref{lin}
  has a unique solution $(u,v)(t,x,y,z)$ satisfying
\[
(u,v)\in L^\infty\Big(0,T;F^{m}_{\alpha,\beta-\rho T}\Big),\quad \pd_z(u,v)\in L^2\Big(0,T;F^{m}_{\alpha,\beta-\rho T}\Big)
\]
with \((u,v)(t,\cdot)\in F^{m}_{\alpha,\beta-\rho t}.
\)
\end{prop}
The proof of this proposition will be given in the Appendix.

From the above Proposition \ref{prop_1} and Proposition \ref{prop_2}, we know that when the monotonic condition is
imposed to one tangential component of background velocity field, such as $u^s$, the analyticity requirement for the velocity field with respect to the corresponding horizontal variable $x$ can be replaced by the
Sobolev regularity, while the velocity field is still analytic in the other horizontal variable $y$, then one still has the local in time well-posedness for the
linearized system in three dimensional space.

Following this argument, it is natural and interesting to study whether the well-posedness of the linearized Prandtl equations \eqref{lin} still holds in the Sobolev
framework if one imposes monotonicity conditions on both tangential velocity components of background state but without analyticity assumption anymore. The study of this paper  gives a negative answer to the question. In fact,  the following theorem  shows a strong linear instability of three dimensional Prandtl equations around basically shear flow in the Sobolev framework except those
with special structure studied in \cite{LWY}.

To state the result, we need some more notations.
 Denote by the operator $T\in\cl(E_{\alpha,\beta},E_{\alpha,\beta'}):$
\begin{equation}\label{def_T}
T(t,s)\big((u_0,v_0)\big)~:=~(u,v)(t,\cdot),
\end{equation}
where $(u,v)$ is the solution of \eqref{lin} with $(u,v)|_{t=s}=(u_0,v_0)$. Introduce the function spaces for $m,\alpha\geq0,$
\begin{equation}\label{def_h}
\h_\alpha^m~:=~H^m\big(\T^2_{x,y};L_\alpha^2(\R^+_z)\big).
\end{equation}
Since the space $E_{\alpha,\beta}$ is dense in the space $\h_\alpha^m,$ we can
extend the operator $T$ from the space $E_{\alpha,\beta}$ 
to $\h_\alpha^m$, 
and define
\[\|T(t,s)\|_{\cl(\h_\alpha^{m_1},\h_\alpha^{m_2})}~:=~\sup\limits_{(u_0,v_0)\in E_{\alpha,\beta}}\frac{\|T(t,s)(u_0,v_0)\|_{\h_\alpha^{m_2}}}{\|(u_0,v_0)\|_{\h^{m_1}_\alpha}}
\in~\R^+\cup\{\infty\},\]
where the infinity means that $T$ can not be extended to $\cl(\h_\alpha^{m_1},\h_\alpha^{m_2})$.

The main result on linear instability of shear flow is stated as follows.

\begin{theorem}\label{thm_lin}
Let $(u^s,v^s)(t,z)$ be the solution of the problems \eqref{shear} satisfying
 $$(u^s-U_0,v^s-V_0)\in C^0\big(\R^+; W_\alpha^{4,\infty}(\R^+)\cap H_\alpha^4(\R^+)\big)\cap C^1\big(\R^+; W_\alpha^{2,\infty}(\R^+)\cap H_\alpha^2(\R^+)\big).$$
Assume that the initial data of \eqref{shear} satisfies that
\begin{equation}\label{ass}
\exists~ z_0>0,~s.t.~(U_s'(z_0))^2+(V_s'(z_0))^2\neq0,~V_s'(z_0)U_s''(z_0)\neq U_s'(z_0)V_s''(z_0).
\end{equation}
Then we have the following two instability statements.

i). There exists $\sigma>0$ such that for all $\delta>0$,
\begin{equation}\label{est_in}
\sup\limits_{0\leq s\leq t\leq\delta}\big\|e^{-\sigma(t-s)\sqrt{|\pd_\cT|}}T(t,s)\big\|_{\cl(\h_\alpha^m,\h_\alpha^{m-\mu})}
~=~+\infty,\quad \forall m>0,~\mu\in[0,\frac{1}{4}),
\end{equation}
where the operator $\pd_\cT$ represents the tangential derivative $\pd_x$ or $\pd_y$;

ii).  There exists an initial shear layer $(U_s,V_s)$ to \eqref{shear} and $\sigma>0$, such that for all $\delta>0$,
\begin{equation}\label{est_in1}
\sup\limits_{0\leq s\leq t\leq\delta}\big\|e^{-\sigma(t-s)\sqrt{|\pd_\cT|}}T(t,s)\big\|_{\cl(\h_\alpha^{m_1},\h_\alpha^{m_2})}
~=~+\infty,\quad \forall m_1,m_2>0.
\end{equation}
\end{theorem}

\begin{remark}\label{rem_3d}
From Theorem \ref{thm_lin}, we know that the three dimensional Prandtl equations can be linearly unstable around the shear flow $(u^s,v^s,0)(t,z)$ even under the monotonic conditions $u_z^s>0$ and $v_z^s>0$.
On the other hand, if we impose
the monotonic condition $U_s'>0$, then \eqref{ass} is equivalent to
\begin{equation}\label{ass1}
\frac{d}{dz}\big(\frac{V_s'}{U_s'}\big)\not\equiv 0.
\end{equation}
And then, by virtue of the boundary condition $U_s(0)=V_s(0)=0$, \eqref{ass1} is equivalent to
\begin{equation}
\frac{d}{dz}\big(\frac{V_s}{U_s}\big)\not\equiv 0.
\end{equation}
Thus,  the result of Theorem \ref{thm_lin} is exactly a complement to the well-posedness result of the
three dimensional Prandtl equations obtained by the authors  in \cite{LWY} for
flow with special structure, that is $\frac{d}{dz}\big(\frac{V_s}{U_s}\big)\equiv 0$. For simplicity, we will  assume that $U_s'(z_0)\neq0$ in the following argument.
\end{remark}

Finally, under the above assumption \eqref{ass} we shall have nonlinear instability for the  original  problem \eqref{3dpd} of the three dimensional nonlinear Prandtl equations.
To state the result, let us first recall the definition of local
 well-posedness from  \cite{guo}.

\begin{define}\label{def_non}
The problem \eqref{3dpd} with the initial data $(u,v)|_{t=0}=(u_0,v_0)(x,y,z)$ is locally well-posed, if there exist positive continuous functions $T(\cdot,\cdot),C(\cdot,\cdot)$, some $\alpha>0$ and some integer $m\geq1$ such that for any initial data $(u_0^1,v_0^1)$ and $(u_0^2,v_0^2)$ with
\[(u_0^1-U_0,v_0^1-V_0)\in \h_\alpha^m,\qquad (u_0^2-U_0,v_0^2-V_0)\in \h_\alpha^m,\]
there are unique distributional solutions $(u^1,v^1)$ and $(u^2,v^2)$ satisfying that for $i=1,2,$
\((u^i,v^i)|_{t=0}=(u_0^i,v_0^i)\)
 and
 \[(u^i-U_0,v^i-V_0)\in L^\infty\big(0,T;L^2(\T^2\times\R^+)\big)\cap L^2\big(0,T;H^1(\T^2\times\R^+)\big),\quad i=1,2,\]
and the following estimate holds
\begin{equation}\label{non}\begin{split}
&\|(u^1,v^1)-(u^2,v^2)\|_{L^\infty(0,T;L^2(\T^2\times\R^+))}+\|(u^1,v^1)-(u^2,v^2)\|_{L^2(0,T;H^1(\T^2\times\R^+))}\\
&\leq C\Big(\|(u_0^1-U_0,v_0^1-V_0)\|_{\h_\alpha^m},\|(u_0^1-U_0,v_0^1-V_0)\|_{\h_\alpha^m}\Big)\|(u_0^1-u_0^2,v_0^1-v_0^2)\|_{\h_\alpha^m},
\end{split}\end{equation}
where $T=T\Big(\|(u_0^1-U_0,v_0^1-V_0)\|_{\h_\alpha^m},\|(u_0^1-U_0,v_0^1-V_0)\|_{\h_\alpha^m}\Big).$
\end{define}

The second main result of this paper is the following ill-posedness of the nonlinear  problem \eqref{3dpd}.

\begin{theorem}\label{thm_non}
Under the same assumption as given in Theorem \ref{thm_lin}, the problem \eqref{3dpd} of the three dimensional nonlinear Prandtl equations is not locally well-posed in the sense of Definition \ref{def_non}.
\end{theorem}

\section{Linear instability}

In this section, we will prove Theorem \ref{thm_lin} to show the linear instability of three dimensional Prandtl equations. The proof is divided into the following four subsections.

\subsection{The linear instability mechanism}

In this subsection, we develop the method introduced in \cite{GV-D} to  analyse
the linear  instability mechanism of three dimensional Prandtl equations. More precisely,  we will find
some high frequency mode in the tangential variables that grow
exponentially in time.
To illustrate this kind of instability mechanism,
as in \cite{GV-D}, we first replace the background shear flow in \eqref{lin} by its initial data so that the background
profile is independent of time. Corresponding to \eqref{lin}, let us consider the following problem:
\begin{equation}\label{sim_lin}\begin{cases}
\pd_t u+(U_s\pd_x+V_s\pd_y)u+wU_s'-\pd_z^2u=0,\quad&{\rm in}\quad \Omega,\\
\pd_t v+(U_s\pd_x+V_s\pd_y)v+wV_s'-\pd_z^2v=0,\quad&{\rm in}\quad \Omega,\\
\pd_xu+\pd_yv+\pd_zw=0,\quad&{\rm in}\quad \Omega,\\
(u,v,w)|_{z=0}=0,\qquad\limz (u,v)=0.
\end{cases}\end{equation}
Noting that the coefficients in \eqref{sim_lin} are also independent of the tangential variables $x$ and $y$, it is convenient to work on
 the Fourier variables with respect to $x$ and $y$.

Recalling in Remark \ref{rem_3d} we assume that $U_s'(z_0)\neq 0$, so we can set \begin{equation}\label{def_a}
a\triangleq\frac{V_s'(z_0)}{U_s'(z_0)},
\end{equation}
and then, the condition
\eqref{ass} yields that the initial tangential velocity $V_s-a U_s$ has a non-degenerate critical point at $z=z_0$. Thus, we may look for solutions of \eqref{sim_lin} in the form of
\begin{equation}\label{expr}
(u,v,w)(t,x,y,z)~=~e^{ik(y-ax+\lambda(k)t)}(\hat u^k,\hat v^k,\hat w^k)(z),
\end{equation}
with some large integer $k$. To insure that the right hand side of \eqref{expr} is $2\pi-$periodic both in $x$ and $y$, that is, both of $k$ and $ak$ being integers, we need that the constant $a$ given in \eqref{def_a} is a rational number, this condition can be easily satisfied. Indeed, the assumption $U_s'(z_0)\neq 0$ and condition \eqref{ass} imply that
\begin{equation}
~U_s'(z)\neq0,\quad V_s'(z)U_s''(z)\neq U_s'(z)V_s''(z)
\end{equation}
holds in a neighborhood of $z_0$. So, by using the  continuity of $(V_s',U_s', V_s'',U_s'')$ and the denseness of the rational numbers in $\R$, there is a point $z_1$ in the neighborhood of $z_0$, such that the condition \eqref{ass} holds for $z_1$ and $\frac{V_s'(z_1)}{U_s'(z_1)}$ is a rational number.
Therefore, in the following discussion we can always assume that
\begin{equation}\label{def-a}
a~=~
\frac{l}{q},
 \end{equation}
for some co-prime integers $l$ and $q$.

  Letting $\ep\triangleq \frac{1}{k}\ll1$,  combining \eqref{expr} with the divergence free condition in \eqref{sim_lin}, we rewrite \eqref{expr} in the form:
\begin{equation}\label{sol}\begin{cases}
(u,v)(t,x,y,z)~&=~e^{i\epsilon^{-1}(y-ax+\lambda_\ep t)}(u_\ep,v_\ep)(z),\\
\quad w~(t,x,y,z)~&=~-i\ep^{-1}e^{i\ep^{-1}(y-ax+\lambda_\ep t)}w_\ep(z).
\end{cases}\end{equation}
Then, the divergence free condition in \eqref{sim_lin} yields that
\begin{equation}\label{diver}
-au_\ep+v_\ep~=~w_\ep'.
\end{equation}
Substituting \eqref{sol} into \eqref{sim_lin}, we obtain
\begin{equation}\label{fr_lin}\begin{cases}
(\lambda_\ep+V_s-aU_s)u_\ep-U_s'w_\ep+i\ep u_\ep^{(2)}=0,\\
(\lambda_\ep+V_s-aU_s)v_\ep-V_s'w_\ep+i\ep v_\ep^{(2)}=0,\\
(u_\ep,v_\ep,w_\ep)(0)=0,\quad\limz (u_\ep,v_\ep)=0.
\end{cases}\end{equation}
Set
\begin{equation}\label{def_wa}
W_s(z)~\triangleq~ (V_s-aU_s)(z),
\end{equation}
combining \eqref{diver} with \eqref{fr_lin} implies that
\begin{equation}\label{fr_w}\begin{cases}
(\lambda_\ep+W_s)w_\ep'-W_s'w_\ep+i\ep w_\ep^{(3)}=0,\\
w_\ep(0)=w_\ep'(0)=0.
\end{cases}\end{equation}
Note that the equation for $w_\ep(z)$ in \eqref{fr_w} is the same as (2.3) studied in \cite{GV-D}
for two dimensional Prandtl equations. Therefore, according to \cite{GV-D}, we have the following result:
\begin{lemma}\label{lem_2d}
For the equation \eqref{fr_w}, if $z_0$ is a non-degenerate critical point of $W_s(z)$, then $w_\ep(z)$ has the following formal approximate expansion in $\ep:$
\begin{equation}\label{ex_w}\begin{cases}
\lambda_\ep~&\sim~-W_s(z_0)+\ep^{\frac{1}{2}}\tau,\\
w_\ep(z)~&\sim~
H(z-z_0)\big[W_s(z)-W_s(z_0)+\ep^{\frac{1}{2}}\tau\big]+\ep^{\frac{1}{2}}W(\frac{z-z_0}{\ep^{\frac{1}{4}}}),
\end{cases}\end{equation}
where $H$ is the Heaviside function, $\tau$ is a complex constant with $\Im\tau<0$, and
 the function $W(Z)$ solves the following ODE:
\begin{equation}\label{eq_w}\begin{cases}
\Big(\tau+W_s''(z_0)\frac{Z^2}{2}\Big)W'-W_s''(z_0)ZW+iW^{(3)}=0,\quad Z\neq0,\\
[W]\big|_{Z=0}=-\tau,~[W']\big|_{Z=0}=0,~[W'']\big|_{Z=0}=-W_s''(z_0),\\
\lim\limits_{Z\rightarrow\pm\infty} W~=~0,\quad  exponentially,
\end{cases}\end{equation}
where the notation $[u]\big|_{Z=0}=\lim\limits_{\delta_1\to 0+}u(\delta_1)-\lim\limits_{\delta_2\to 0-}u(\delta_2)$ denotes the jump of a related function $u(Z)$ across $Z=0$.
\end{lemma}

As in \cite{GV-D}, the asymptotic expansion of the solution $w_\ep(z)$ given in \eqref{ex_w} shows that
the approximate solution of $w_\ep(z)$ can be divided into the "regular" part and the "shear layer" part as:
\[ w_\ep^a(z)~:=~w_\ep^{reg}(z)+w_\ep^{sl(z)}\]
with
\[w_\ep^{reg}~\triangleq~H(z-z_0)\big[W_s(z)-W_s(z_0)+\ep^{\frac{1}{2}}\tau\big],
\]
and
\[w_\ep^{sl}~\triangleq~\ep^{\frac{1}{2}}W(\frac{z-z_0}{\ep^{\frac{1}{4}}}).\]
Note that $w_\ep^{reg}(0)=(w_\ep^{reg})'(0)=0.$ The "shear layer" part $w_\ep^{sl}$ is to cancel  the discontinuities of the "regular" part $w_\ep^{reg}$ at $z=z_0$, such that the approximation $w_\ep^{a}\in C^2(\R^+)$.

The formal asymptotic expansion $\eqref{ex_w}_1$ for the eigenvalue indicates strong instability of \eqref{sim_lin}, that is, back to the Fourier representation \eqref{sol}, the tangential velocity $(u,v)$ grows like $e^{\frac{t}{\sqrt{\ep}}}$. To complete this process, we will construct the formal approximation of $(u_\ep,v_\ep)(z)$.
The construction of $(u_\ep,v_\ep)(z)$ is based on the relation \eqref{diver} and the approximation \eqref{ex_w}, which implies that
\begin{equation}\label{ex-w}
(v_\ep-au_\ep)(z)~\sim~H(z-z_0)W_s'(z)+\ep^{\frac{1}{4}}W'(\frac{z-z_0}{\ep^{\frac{1}{4}}}).
\end{equation}
From \eqref{ex-w}, we assume that
the formal approximation for $(u_\ep,v_\ep)(z)$ can be chosen as follows:
\begin{equation}\label{ex_uv}\begin{cases}
u_\ep(z)~\sim~H(z-z_0)U_s'(z)+h(\frac{z-z_0}{\ep^{\frac{1}{4}}})+\ep^{\frac{1}{4}}U(\frac{z-z_0}{\ep^{\frac{1}{4}}}),\\
v_\ep(z)~\sim~H(z-z_0)V_s'(z)+ah(\frac{z-z_0}{\ep^{\frac{1}{4}}})+\ep^{\frac{1}{4}}\big(aU+W'\big)(\frac{z-z_0}{\ep^{\frac{1}{4}}}),
\end{cases}\end{equation}
where the functions $h(Z)$ and $U(Z)$ is rapidly decay as $Z\rightarrow\pm\infty$
as explained in the following.
One can easily verifies that
\[H(z-z_0)\big(U_s'(z),V_s'(z)\big)\big|_{z=0}=0,\qquad\limz (u_\ep,v_\ep)~=~0.\]
Similar to the ``shear layer'' part $w_\ep^{sl}$ defined
in \eqref{ex_w}, the functions $h(Z)$ and $U(Z)$ are used
to cancel the discontinuities in $H(z-z_0)U_s'(z)$ and $H(z-z_0)V_s'(z)$ at $z=z_0$,
so that the approximations in \eqref{ex_uv} belongs to $C^1(\R^+)$. Moreover, $h(Z)$
and $U(Z)$ are used to balance the approximation in the orders
of  $O(\sqrt{\ep})$ and $O(\ep^{\frac{3}{4}})$ respectively.
For this, from \eqref{fr_lin}, \eqref{ex_w} and \eqref{ex_uv},  $h(Z)$ and $U(Z)$ satisfy the following problems respectively,
\begin{equation}\label{eq_h}\begin{cases}
\Big(\tau+W_s''(z_0)\frac{Z^2}{2}\Big)h-U_s'(z_0)W+ih''=0,\quad Z\neq0,\\
[h]\big|_{Z=0}=-U_s'(z_0),\quad[h']\big|_{Z=0}=0,\\
\limZ h~=~0,
\end{cases}\end{equation}
and
\begin{equation}\label{eq_u}\begin{cases}
\Big(\tau+W_s''(z_0)\frac{Z^2}{2}\Big)U-U_s''(z_0)ZW+iU''=0,\quad Z\neq0,\\
[U]\big|_{Z=0}=0,\qquad [U']\big|_{Z=0}=-U_s''(z_0),\\
\limZ U~=~0.
\end{cases}\end{equation}

Comparing the problem \eqref{eq_h} with \eqref{eq_u}, respectively \eqref{eq_w}, it follows that $\frac{U_s'(z_0)}{W_s''(z_0)}W''(Z)$, respectively $\frac{U_s''(z_0)}{W_s''(z_0)}W'(Z)$, solves the problem \eqref{eq_h}, respectively \eqref{eq_u}. Consequently,
we can choose the formal expansion of $(u_\ep,v_\ep)(z)$ as
\begin{equation}\label{ex_uv1}\begin{cases}
u_\ep(z)~\sim~H(z-z_0)U_s'(z)+\frac{U_s'(z_0)}{W_s''(z_0)}W''(\frac{z-z_0}{\ep^{\frac{1}{4}}})+\ep^{\frac{1}{4}}\frac{U_s''(z_0)}{W_s''(z_0)}W'(\frac{z-z_0}{\ep^{\frac{1}{4}}}),\\
v_\ep(z)~\sim~H(z-z_0)V_s'(z)+\frac{V_s'(z_0)}{W_s''(z_0)}W''(\frac{z-z_0}{\ep^{\frac{1}{4}}})+\ep^{\frac{1}{4}}\frac{V_s''(z_0)}{W_s''(z_0)}W'(\frac{z-z_0}{\ep^{\frac{1}{4}}}).
\end{cases}\end{equation}

Therefore, we have concluded the following results for the reduced boundary value problem \eqref{sim_lin}.

\begin{prop}\label{sim-lin}

For the large frequency $k=\frac{1}{\epsilon}$, the approximate solutions of the problem \eqref{sim_lin} can be expressed as \eqref{sol} with
\begin{equation}\label{ex_all}\begin{cases}
\lambda_\ep~&\sim~-W_s(z_0)+\ep^{\frac{1}{2}}\tau,\\
w_\ep(z)~&\sim~H(z-z_0)\Big[W_s(z)-W_s(z_0)+\ep^{\frac{1}{2}}\tau\Big]
+\ep^{\frac{1}{2}}W(\frac{z-z_0}{\ep^{\frac{1}{4}}}),\\
(u_\ep(z),v_\ep(z))~&\sim~H(z-z_0)(U_s'(z),V_s'(z))+
W''(\frac{z-z_0}{\ep^{\frac{1}{4}}})\frac{1}{W_s''(z_0)}\big(U_s'(z_0),V_s'(z_0)\big),\\
&\quad+\ep^{\frac{1}{4}}W'(\frac{z-z_0}{\ep^{\frac{1}{4}}})\frac{1}{W_s''(z_0)}\big(U_s''(z_0),V_s''(z_0)\big),
\end{cases}\end{equation}
where the complex constant $\tau$ and function $W(Z)$ are given in Lemma \ref{lem_2d}.

\end{prop}

\begin{remark}\label{rem_tau}

Recalling from \cite{GV-D}, we know that the pair $\big(\tau, W(Z)\big)$ given in Lemma \ref{lem_2d} has the following form
\begin{equation}\label{ch_w}\begin{cases}
\tau~&=~\big|\frac{W_s''(z_0)}{2}\big|^{\frac{1}{2}}~\tilde\tau,\\
W(Z)~&=~\big|\frac{W_s''(z_0)}{2}\big|^{\frac{1}{2}}\Big[\big(\tilde\tau+\big|\frac{W_s''(z_0)}{2}\big|^{\frac{1}{2}}Z^2\big)\tilde W\big(\big|\frac{W_s''(z_0)}{2}\big|^{\frac{1}{4}}Z\big)-1_{\R^+}\big(\tilde\tau+\big|\frac{W_s''(z_0)}{2}\big|^{\frac{1}{2}}Z^2\big)\Big].
\end{cases}\end{equation}
where the complex constant $\tilde\tau$ has a negative imaginary part, i.e., $\Im\tilde\tau<0,$ and the function $\tilde W(\tilde Z)$ is a smooth solution of the following third order ordinary differential equation:
\begin{equation}\label{SC}\begin{cases}
\big(\tilde\tau+sign(W_s''(z_0))\tilde Z^2\big)^2\frac{d}{d\tilde Z}\tilde W+i\frac{d^3}{d\tilde Z^3}\Big(\big(\tau+sign(W_s''(z_0))\tilde Z^2\big)\tilde W\Big)=0,\\
\lim\limits_{\tilde Z\rightarrow-\infty}\tilde W~=~0,\quad\lim\limits_{\tilde Z\rightarrow+\infty}\tilde W~=~1.
\end{cases}\end{equation}

\end{remark}

\subsection{Construction of approximate solutions}

Inspiring by the construction of approximate solutions to the simplified problem \eqref{sim_lin} given in
the above subsection, and also by the argument given in \cite{GV-D}, we are going to construct the approximate solution of the original linearized problem \eqref{lin}.

Let $(u^s(t,z),v^s(t,z))$ satisfy the assumptions of Theorem \ref{thm_lin},
and denote by
\[ w_a^s(t,z)~\triangleq~v^s(t,z)-a u^s(t,z),\]
where the constant $a$ is given in \eqref{def_a}.
Then, we know that $z_0$ is a non-degenerate critical point of $w_a^s(0,z)$. Without loss of generality, we assume that $\pd_z^2w_a^s(0,z_0)<0$, then the differential equation
\begin{equation}\label{def_alpha}\begin{cases}
\pd_t\pd_zw_a^s\big(t,f(t)\big)+\pd_z^2w_a^s\big(t,f(t)\big)f'(t)=0,\\
f(0)~=~z_0,
\end{cases}\end{equation}
defines a non-degenerate critical point $f(t)$ of $w_a^s(t,\cdot)$ when $0<t<t_0$ for some small $t_0>0$. Note that such $f(t)$ can also be determined by the following equation:
\[\frac{d}{dt}\Big(\frac{v_z^s\big(t,f(t)\big)}{u_z^s\big(t,f(t)\big)}\Big)=0,\quad f(0)=z_0.\]

Since the approximation solution of \eqref{sim_lin} given in Proposition \ref{sim-lin} is obtained with the background state being frozen at the initial data $(u^s,v^s)|_{t=0}=(U_s,V_s)(z)$, to construct approximate solutions of the original problem \eqref{lin} with background state being the shear flow in the time interval $0<t<t_0$, we need to do some modification as in \cite{GV-D}. Recall Proposition \ref{sim-lin} and Remark \ref{rem_tau} in the above subsection, let $\tau, W(Z)$ be given in \eqref{ch_w}, and set
\begin{equation}\begin{split}\label{def_phi}
 W_{sl}(Z)~&:=~\Big(\tau-Z^2\Big) W( Z)-1_{\R^+}\Big(\tau-Z^2\Big).
\end{split}\end{equation}
Then, for $\ep>0$ we introduce
\begin{equation}\label{app_lam}
\lambda_\ep(t)~:=~-w_a^s(t,f(t))+\ep^{\frac{1}{2}}\big|\frac{\pd_z^2w_a^s(t,f(t))}{2}\big|^{\frac{1}{2}}~\tau,
\end{equation}
and the ``regular'' part of velocities
\begin{equation}\label{app_reg}\begin{cases}
U_\ep^{reg}(t,z)=H(z-f(t))u_z^s(t,z),\quad
V_\ep^{reg}(t,z)=H(z-f(t))v_z^s(t,z),\\[3mm]
W_\ep^{reg}(t,z)=H(z-f(t))\Big[w_a^s(t,z)-w_a^s(t,f(t))+\ep^{\frac{1}{2}}\big|\frac{\pd_z^2w_a^s(t,f(t))}{2}\big|^{\frac{1}{2}} \tau\Big],
\end{cases}\end{equation}
as well as the ``shear layer'' part of $w_\ep$
\begin{equation}\label{app_sl}\begin{split}
W_\ep^{sl}(t,z)~&:=~\ep^{\frac{1}{2}}\varphi(z-f(t))\big|\frac{\pd_z^2w_a^s(t,f(t))}{2}\big|^{\frac{1}{2}} W_{sl}\Big(\big|\frac{\pd_z^2w_a^s(t,f(t))}{2}\big|^{\frac{1}{4}}\cdot\frac{z-f(t)}{\ep^{\frac{1}{4}}}\Big).
\end{split}\end{equation}
Here, $\varphi$ is a smooth truncation function near 0, and $W_{sl}$ is given in \eqref{def_phi}. Therefore, from Proposition \ref{sim-lin} and by \eqref{app_reg} and \eqref{app_sl}, the approximate solution of the problem \eqref{lin} can be defined as:
\begin{equation}\label{ap_lin}
(u_\ep,v_\ep,w_\ep)(t,x,y,z)=e^{i\ep^{-1}(y-ax)}\big(U_\ep,V_\ep,W_\ep\big)(t,z)
\end{equation}
with
\begin{equation}\begin{split}\label{expre_uv}
\big(U_\ep,V_\ep\big)(t,z)~&=~ie^{i\ep^{-1}\int_0^t\lambda_\ep(s)ds}\Big\{\big(U_\ep^{reg},V_\ep^{reg}\big)(t,z)
+\frac{\pd_z^2W_\ep^{sl}(t,z)}{\pd_z^2w_a^s(t,f(t))}\big(u_z^s,v_z^s\big)(t,f(t))\\ 
&\qquad\qquad\qquad\qquad+\frac{\pd_zW_\ep^{sl}(t,z)}{\pd_z^2w_a^s(t,f(t))}\big(\pd_z^2u^s,\pd_z^2v^s\big)(t,f(t))\Big\},\\
W_\ep(t,z)~&=~\ep^{-1}e^{i\ep^{-1}\int_0^t\lambda_\ep(s)ds}\Big(W_\ep^{reg}(t,z)+W_\ep^{sl}(t,z)\Big).
\end{split}\end{equation}
 Moreover, in order that the function \eqref{ap_lin} is $2\pi-$periodic in $x$ and $y$, we take $\ep=\frac{1}{qk}$ with the integers $q$ given in \eqref{def-a} and $k\in\N$.

It is straightforward to check that for $(u_\ep,v_\ep,w_\ep)$ defined in \eqref{ap_lin},
\[(u_\ep,v_\ep,w_\ep)|_{z=0}=0,\qquad\limz (u_\ep,v_\ep)=0,\]
 and the divergence free condition holds. Also, $$(u_\ep,v_\ep)(t,x,y,z)=e^{i\ep^{-1}(y-ax)}(U_\ep,V_\ep)(t,z)$$ is analytic in the tangential variables $x,y$ and $H_\alpha^{2}$ in $z$. Moreover, there are positive constants $C_0$ and $\sigma_0$, independent of $\ep$, such that
\begin{equation}\label{bound_app}
C_0^{-1}e^{\frac{\sigma_0t}{\sqrt{\ep}}}~\leq~\|(U_\ep,V_\ep)(t,\cdot)\|_{L_\alpha^{2}}
~\leq~C_0e^{\frac{\sigma_0t}{\sqrt{\ep}}}.
\end{equation}

Plugging the relation \eqref{ap_lin} into the original linearized Prandtl equations \eqref{lin}, it follows that
\begin{equation}\label{lin_ap}\begin{cases}
\pd_t u_\ep+(u^s\pd_x+v^s\pd_y)u_\ep+w_\ep u_z^s-\pd_z^2u_\ep=r_\ep^{1},\\
\pd_t v_\ep+(u^s\pd_x+v^s\pd_y)v_\ep+w_\ep v_z^s-\pd_z^2v_\ep=r_\ep^2,\\
\pd_xu_\ep+\pd_yv_\ep+\pd_zw_\ep=0,\quad{\rm in}\quad \Omega,\\
(u_\ep,v_\ep,w_\ep)|_{z=0}=0,\qquad\limz (u_\ep,v_\ep)=0,
\end{cases}\end{equation}
where the remainder term is represented by
\begin{equation}\label{def_err}
(r_\ep^1,r_\ep^2)(t,x,y,z):=~e^{i\ep^{-1}(y-ax)}(R_\ep^1,R_\ep^2)(t,z)
\end{equation}
with
\begin{equation}\label{r}\begin{split}
R_\ep^1&=\pd_t U_\ep+i\ep^{-1}w_a^s(t,z)U_\ep-\pd_z^2U_\ep+u_z^s(t,z)W_\ep,\\
R_\ep^2&=\pd_t V_\ep+i\ep^{-1}w_a^s(t,z)V_\ep-\pd_z^2V_\ep+v_z^s(t,z)W_\ep.
\end{split}\end{equation}
Note that the representation \eqref{app_reg}  implies
\begin{equation}\label{eq_reg}
(\pd_t-\pd_z^2)U_\ep^{reg}=(\pd_t-\pd_z^2)V_\ep^{reg}=0,\quad z\neq f(t).
\end{equation}
Then, from the representation of $(U_\ep,V_\ep)$ and $W_\ep$ given in \eqref{expre_uv} and by using the equations \eqref{SC} and \eqref{eq_reg}, we conclude that for $z\neq f(t),$
\begin{equation}\label{r1}\begin{split}
R_\ep^1(t,z)
&=e^{i\ep^{-1}\int_0^t\lambda_\ep(s)ds}\Big\{-\ep^{-1}\Big[w_a^s(t,z)
-w_a^s(t,f(t))-\pd_z^2w_a^s(t,f(t))\frac{(z-f(t))^2}{2}\Big]\\
&\qquad\qquad\cdot\Big[
\frac{u_z^s(t,f(t))}{\pd_z^2w_a^s(t,f(t))}\pd_z^2W_\ep^{sl}+\frac{\pd_z^2u^s(t,f(t))}{\pd_z^2w_a^s(t,f(t))}\pd_zW_\ep^{sl}\Big]\\
&\qquad+\ep^{-1}\Big[u_z^s(t,z)-u_z^s(t,f(t))-u_{zz}^s(t,f(t))(z-f(t))\Big]W_\ep^{sl}\\
&\qquad+i\pd_t\Big(
\frac{u_z^s(t,f(t))}{\pd_z^2w_a^s(t,f(t))}\pd_z^2W_\ep^{sl}+\frac{\pd_z^2u^s(t,f(t))}{\pd_z^2w_a^s(t,f(t))}\pd_zW_\ep^{sl}\Big)
+O(\ep^\infty)\Big\},
\end{split}\end{equation}
and
\begin{equation}\label{r2}\begin{split}
R_\ep^2(t,z)
&=e^{i\ep^{-1}\int_0^t\lambda_\ep(s)ds}\Big\{-\ep^{-1}\Big[w_a^s(t,z)
-w_a^s(t,f(t))-\pd_z^2w_a^s(t,f(t))\frac{(z-f(t))^2}{2}\Big]\\
&\qquad\qquad\cdot\Big[
\frac{v_z^s(t,f(t))}{\pd_z^2w_a^s(t,f(t))}\pd_z^2W_\ep^{sl}+\frac{\pd_z^2v^s(t,f(t))}{\pd_z^2w_a^s(t,f(t))}\pd_zW_\ep^{sl}\Big]\\
&\qquad+\ep^{-1}\Big[v_z^s(t,z)-v_z^s(t,f(t))-v_{zz}^s(t,f(t))(z-f(t))\Big]W_\ep^{sl}\\
&\qquad+i\pd_t\Big(
\frac{v_z^s(t,f(t))}{\pd_z^2w_a^s(t,f(t))}\pd_z^2W_\ep^{sl}+\frac{\pd_z^2v^s(t,f(t))}{\pd_z^2w_a^s(t,f(t))}\pd_zW_\ep^{sl}\Big)
+O(\ep^\infty)\Big\}.
\end{split}\end{equation}
The terms $O(\ep^\infty)$ in \eqref{r1} and \eqref{r2} represent the remainders
with exponential decay in $z$ that follows from the fact that $W_\ep^{sl}$ decays exponentially and the derivatives of $\varphi(\cdot-f(t))$ vanish near $f(t)$.  Then,  with the same $\sigma_0$ given in \eqref{bound_app}, we have that
$(R_\ep^1,R_\ep^2)(t,z)$ satisfy
\begin{equation}\label{bound_r}
\|(R_\ep^1,R_\ep^2)(t,\cdot)\|_{L_\alpha^{2}}\leq C_1\ep^{-\frac{1}{4}}e^{\frac{\sigma_0t}{\sqrt{\ep}}},
\end{equation}
where the constant $C_1$ is independent of $\ep$.

Therefore, we conclude
\begin{prop}\label{con_app}
For the linear problem \eqref{lin}, the approximate solution $(u_\ep,v_\ep,w_\ep)$ given in \eqref{ap_lin} and \eqref{expre_uv}, satisfies the problem \eqref{lin_ap}, where the source term $(r_\ep^1,r_\ep^2)$ in the form \eqref{def_err} has the bound \eqref{bound_r}.
\end{prop}

\begin{remark}\label{rem_r}
The estimate \eqref{bound_r} follows from the fact that the "shear layer" part $\pd_z^2W_\ep^{sl}$ cancels the terms $\ep^{-1}u_z^s(t,f(t))W_\ep^{sl}$ and $\ep^{-1}v_z^s(t,f(t))W_\ep^{sl}$ in \eqref{r1} and \eqref{r2} respectively, by using the equation \eqref{eq_h} for instance. And this
error bound leads to the choice of  $\mu<\frac{1}{4}$ in \eqref{est_in}. It is slight different from the two dimensional case studied in \cite{GV-D} where  $\mu<\frac{1}{2}$, because here we require that the initial data of $u^s$ and $v^s$ do not degenerate simultaneously at a point, $z=z_0$, while in the two dimensional problem, it is assumed that there is a degeneracy at the critical point.
\end{remark}

\subsection{Proof of Theorem \ref{thm_lin}(i)}

At this stage, based on the approximate solution given in Proposition \ref{con_app}, we can use the method from \cite{GV-D} to prove  Theorem \ref{thm_lin}. We
now  sketch the proof as follows.

\underline{\it Verification of \eqref{est_in} for the tangential differential operator by contradiction}.
Suppose that \eqref{est_in} does not hold for $\pd_x$, that is, for all $\sigma>0$, there exists $\delta>0,m\geq0$ and $\mu\in[0,\frac{1}{4})$, that is,
\begin{equation}\label{pr_x}
\sup\limits_{0\leq s\leq t\leq \delta}\|e^{-\sigma(t-s)\sqrt{|\pd_x|}}T(t,x)\|_{\cl(\h_\alpha^m,\h_\alpha^{m-\mu})}<+\infty.
\end{equation}
Introduce the operator
\begin{equation*}
T_\ep(t,s):~L_{\alpha}^2(\R^+)\mapsto L_{\alpha}^2(\R^+)
\end{equation*}
as
\begin{equation}\label{def_t}
 T_\ep(t,s)\big((U_0,V_0)\big)~:=~e^{-i\ep^{-1}(y-ax)}T(t,s)\Big(e^{i\ep^{-1}(y-ax)}(U_0,V_0)\Big)
\end{equation}
with $T(t,s)$ being defined in \eqref{def_T}.
From \eqref{pr_x}, we have
\begin{equation}\label{est_tep}
\|T_\ep(t,s)\|_{\cl(L_\alpha^{2})}~\leq~C_2\ep^{-\mu}e^{\frac{\sqrt{a}\sigma(t-s)}{\sqrt{\ep}}},\quad \forall ~0\leq s\leq t\leq\delta
\end{equation}
for a constant $C_2$ independent of $\ep$.

Denote by
\[L_\ep~:=~e^{-i\ep^{-1}(y-ax)}~L~e^{i\ep^{-1}(y-ax)},\]
where $L$ is the linearized Prandtl operator around the shear flow $(u^s,v^s,0)$.
Let $(U,V)(t,z)$ be a solution to the problem
\[\begin{cases}
\pd_t(U,V)+L_\ep(U,V)~=~0,\\
(U,V)|_{t=0}~=~(U_\ep,V_\ep)(0,z).
\end{cases}\]

From the definition \eqref{def_t}, we have that
\[(U,V)(t,z)~=~T_\ep(t,0)\big((U_\ep,V_\ep)(0,z)\big).\]
Then, from \eqref{est_tep} it follows that for all $t\leq\delta,$
\begin{equation}\label{up_bound}
\|(U,V)(t,\cdot)\|_{L_\alpha^{2}}\leq C_2\ep^{-\mu}e^{\frac{\sqrt{a}\sigma t}{\sqrt{\ep}}}\|(U_\ep,V_\ep)(0,\cdot)\|_{L_\alpha^{2}}
\leq C_3\ep^{-\mu}e^{\frac{\sqrt{a}\sigma t}{\sqrt{\ep}}}
\end{equation}
holds for a constant $C_3$ independent of $\ep$. From \eqref{r1} and \eqref{r2},
we know that
\[\pd_t(U_\ep,V_\ep)+L_\ep(U_\ep,V_\ep)~=~(R_\ep^1,R_\ep^2).\]
Thus, the difference $(\tilde U,\tilde V):=(U,V)-(U_\ep,V_\ep)$ can be obtained by the Duhamel representation:
\begin{equation}\label{eq_dif}
(\tilde U,\tilde V)(t,\cdot)~=~\int_0^tT_\ep(t,s)\big((R^1_\ep,R_\ep^2)(s,\cdot)\big)ds,\quad \forall~t\leq\delta.
\end{equation}
Combining \eqref{bound_r}, \eqref{est_tep} and \eqref{eq_dif},
and choosing $\sqrt{a}\sigma<\sigma_0$ yields that
\begin{equation}\label{est_dif}
\|(\tilde U,\tilde V)(t,\cdot)\|_{L^2_\alpha}\leq C_1C_2\ep^{-\mu-\frac{1}{4}}
\int_0^te^{\frac{\sqrt{a}\sigma(t-s)}{\sqrt{\ep}}}e^{\frac{\sigma_0s}{\sqrt{\ep}}}ds
\leq C_4\ep^{\frac{1}{4}-\mu}e^{\frac{\sigma_0t}{\sqrt{\ep}}},
\end{equation}
where the constant $C_4>0$ is independent of $\ep$. Then, by using \eqref{bound_app},
we obtain that for $t<\delta$ and sufficiently small $\ep$,
\begin{equation}\label{low_bound}\begin{split}
\|(U,V)(t,\cdot)\|_{L_\alpha^2}&\geq\|(U_\ep,V_\ep)(t,\cdot)\|_{L_\alpha^2}
-\|(\tilde U,\tilde V)(t,\cdot)\|_{L_\alpha^2}\\
&\geq C_0^{-1}e^{\frac{\sigma_0t}{\sqrt{\ep}}}-C_4\ep^{\frac{1}{4}-\mu}e^{\frac{\sigma_0t}{\sqrt{\ep}}}
\geq \frac{1}{2}C_0^{-1}e^{\frac{\sigma_0t}{\sqrt{\ep}}}.
\end{split}\end{equation}
As $\sqrt{a}\sigma<\sigma_0,$ comparing \eqref{up_bound} with \eqref{low_bound},
the contradiction arises when $t>\frac{\sqrt{\ep}}{\sigma_0-\sqrt{a}\sigma}\big(\ln(2C_0C_3)-\ln\ep\big)$
with sufficiently small $\ep$. Thus, the proof of Theorem \ref{thm_lin}(i) is completed.

\subsection{Proof of Theorem \ref{thm_lin}(ii)}

The aim of this subsection is to prove Theorem \ref{thm_lin}(ii). By comparing \eqref{est_in} with \eqref{est_in1}, we only need to show that there exists an initial data $(U_s,V_s)$ for
the shear flow to \eqref{shear} such that \eqref{est_in} still holds for arbitrary $\mu>0.$
Recall the  proof of part i) in the above subsection, the task can be attributed to find some $(U_s,V_s)$ such that the remainder $(R_\ep^1,R_\ep^2)(t,z)$, generated in \eqref{r} by the approximation \eqref{expre_uv}, has the following estimate:
\begin{equation}\label{est_r}
\|(R_\ep^1,R_\ep^2)(t,\cdot)\|_{L_\alpha^{2}}\leq C\big(\ep^{N}+t^{2N}\big)e^{\frac{\sigma_0t}{\sqrt{\ep}}}
\end{equation}
for some $N\in\R^+,N+\frac{1}{2}>\mu$. Once this is achieved, as in \cite{GV-D}, the desired conclusion holds.

Similar to \cite{GV-D}, the special initial shear layer $(U_s,V_s)(z)$ can be chosen such that
\[\begin{split}
(U_s,V_s)|_{z=0}=0,\quad(U_s,V_s)(z)\rightarrow(U_0,V_0),~{\rm exponentially~as}~z\rightarrow+\infty,
\end{split}\]
and in a small neighborhood of $z_0$,
\begin{equation}\label{sp_uv}\begin{cases}
U_s(z)=U_s''(z_0)\frac{(z-z_0)^2}{2}+q(z-z_0)+U_s(z_0),\\ V_s(z)=V_s''(z_0)\frac{(z-z_0)^2}{2}+l(z-z_0)+V_s(z_0),
\end{cases}\end{equation}
where $q,l$ are integers given in \eqref{def-a}, and the constants $U_s''(z_0),U_s(z_0),$ $V_s''(z_0),V_s(z_0)$ satisfy
\[V_s''(z_0)-aU_s''(z_0)\neq0,\quad U_s(z_0)\neq0\]
with the constant $a$ being given in \eqref{def_a}.
Then, for such $(U_s,V_s)(z)$, we will show that \eqref{est_r} holds.
We only estimate the term
$R_\ep^1$, as the same argument works for $R_\ep^2$. From \eqref{r1},
decompose $R_\ep^1$ as follows for $z\neq f(t)$:
\begin{equation}\label{div}
R_\ep^1~:=~e^{i\ep\int_0^t\lambda(s)ds}\Big(R_{\ep,1}^1+R_{\ep,2}^1+R_{\ep,3}^1+O(\ep^\infty)\Big),
\end{equation}
where
\[\begin{split}
R_{\ep,1}^1(t,z)&=-\ep^{-1}\Big[w_a^s(t,z)
-w_a^s(t,f(t))-\pd_z^2w_a^s(t,f(t))\frac{(z-f(t))^2}{2}\Big]\\
&\qquad\qquad\cdot\Big[
\frac{u_z^s(t,f(t))}{\pd_z^2w_a^s(t,f(t))}\pd_z^2W_\ep^{sl}+\frac{\pd_z^2u^s(t,f(t))}{\pd_z^2w_a^s(t,f(t))}\pd_zW_\ep^{sl}\Big],\\
R_{\ep,2}^1(t,z)&=\ep^{-1}\Big[u_z^s(t,z)-u_z^s(t,f(t))-u_{zz}^s(t,f(t))(z-f(t))\Big]W_\ep^{sl},\\
R_{\ep,3}^1(t,z)&=i\pd_t\Big(\frac{u_z^s(t,f(t))}{\pd_z^2w_a^s(t,f(t))}\pd_z^2W_\ep^{sl}+\frac{\pd_z^2u^s(t,f(t))}{\pd_z^2w_a^s(t,f(t))}\pd_zW_\ep^{sl}\Big).
\end{split}\]
Therefore, it remains to show that
\begin{equation}\label{bound_r1}
\|R_{\ep,i}^1(t,\cdot)\|_{L_\alpha^{2}}~\leq~C(\ep^N+t^{2N}),\quad i=1,2,3.
\end{equation}

Firstly, for $R^1_{\ep,1}$, note that the function $w_a^s(t,z)=v^s(t,z)-au^s(t,z)$ satisfies
\begin{equation}\label{heat}
\pd_t w_a^s-\pd_z^2 w_a^s=0,
\end{equation}
 and in a small neighborhood of $z_0$,
\begin{equation}\label{sp_w}
w_a^s(0,z)=\big(V_s''(z_0)-aU_s''(z_0)\big)\frac{(z-z_0)^2}{2}+V(z_0)-aU_s(z_0).
\end{equation}
Thus, by using the Taylor expansion and \eqref{heat}, it follows that for any $N\in\N$,
\begin{equation}\label{est-r1}\begin{split}
&|R_{\ep,1}^1(t,z)|=\ep^{-1}\Big|\int_{f(t)}^z\frac{(\tilde z-f(t))^2}{2}\pd_z^3w_a^s(t,\tilde z)d\tilde z\Big|\cdot\Big|
\frac{u_z^s(t,f(t))}{\pd_z^2w_a^s(t,f(t))}\pd_z^2W_\ep^{sl}+\frac{\pd_z^2u^s(t,f(t))}{\pd_z^2w_a^s(t,f(t))}\pd_zW_\ep^{sl}\Big|\\
&\leq C\ep^{-1}\int_{f(t)}^z\frac{(\tilde z-f(t))^2}{2}\sum_{k=0}^{2N-1}\frac{t^k}{k!}\big|\pd_t^k\pd_z^3w_a^s(0,\tilde z)\big|d\tilde z\cdot\Big(\big|\pd_z^2W_\ep^{sl}\big|+\big|\pd_zW_\ep^{sl}\big|\Big)+O(t^{2N})\\
&\leq C\ep^{-1}\int_{f(t)}^z\frac{(\tilde z-f(t))^2}{2}\sum_{k=0}^{2N-1}\frac{t^k}{k!}\big|\pd_z^{3+2k}w_a^s(0,\tilde z)\big|d\tilde z\cdot\Big(\big|\pd_z^2W_\ep^{sl}\big|+\big|\pd_zW_\ep^{sl}\big|\Big)+O(t^{2N}).
\end{split}\end{equation}
From \eqref{sp_w}, we know that
\[\pd_z^{3+2k}w_a^s(0,z)~\equiv~0,\quad\quad \forall k\in\N,\]
when $z$ is in a small neighborhood  of $f(t)$ as $t>0$ is small. This implies
 that the integral in the last line of \eqref{est-r1} is supported away from $z=f(t)$.
Then, combining with the exponential decrease of the ``shear layer'' $W_\ep^{sl}$, \eqref{est-r1} yields
\begin{equation}\label{est_r1}
\|R_{\ep,1}^1(t,\cdot)\|_{L_\alpha^{2}}~\leq~C (\ep^N+t^{2N}).
\end{equation}
For the term $R_{\ep,2}^1$, we can use similar arguments as above to obtain
\begin{equation}\label{est_r2}
\|R_{\ep,2}^1(t,\cdot)\|_{L_\alpha^{2}}~\leq~C (\ep^N+t^{2N}).
\end{equation}

Next, from the expression \eqref{app_sl} of $W_\ep^{sl}(t,z)$ and the relation
\[\pd_tu^s=\pd_z^{2}u^s,\quad \pd_tw_a^s=\pd_z^{2}w_a^s,\]
a straightforward calculation implies that
 \begin{equation}\label{expre_r3}\begin{split}
 R_{\ep,3}^1(t,z)=&\pd_z^3u^s(t,f(t))g_1(t,z)+\pd_z^4u^s(t,f(t))g_2(t,z)+
 \pd_z^4w_a^s(t,f(t))g_3(t,z)\\
 &+f'(t)g_4(t,z)+O(\ep^\infty),
\end{split}\end{equation}
where  each function $g_i(t,z),1\leq i\leq4$ can be expressed as a linear combination of the terms
 \[
 O(\ep^{-\frac{1}{4}})\varphi(z-f(t))~W_{sl}^{(j)}\Big(\big|\frac{\pd_z^2w_a^s(t,f(t))}{2}\big|^{\frac{1}{4}}\cdot\frac{z-f(t)}{\ep^{\frac{1}{4}}}\Big),\quad 0\leq j\leq 3
 \]
 with $W_{sl}$ being defined in \eqref{def_phi}. Then, as for $R_{\ep,1}^1$, the first three terms on the right hand side of \eqref{expre_r3} have the same bounds as in \eqref{est_r1}. For the fourth term given in \eqref{expre_r3}, by noticing that from \eqref{def_alpha},
\[
f'(t)~=~-\frac{\pd_t\pd_zw_a^s(t,f(t))}{\pd_z^2w_a^s(t,f(t))}
~=~\frac{\pd_z^3w_a^s(t,f(t))}{\pd_z^2w_a^s(t,f(t))},
\]
we can verify that $f'(t)g_4(t,z)$ also satisfies the same estimate as \eqref{est_r1}. In conclusion, we have
\begin{equation}\label{est_r3}
\|R_{\ep,3}^1(t,\cdot)\|_{L_\alpha^{2}}~\leq~C (\ep^N+t^{2N}).
\end{equation}
Thus, combining \eqref{est_r1}, \eqref{est_r2} and \eqref{est_r3}, we have the estimate \eqref{bound_r1}, and then obtain the proof of Theorem \ref{thm_lin}(ii) by taking $N$ large enough.

Finally, we state the following result to finish this section, which can be obtained by similar arguments as above,
\begin{prop}\label{prop_uni}
There exists a shear layer $(u^s,v^s)$ to \eqref{shear} with $u_z^s,v_z^s>0$, such that for all $\delta>0$,
\begin{equation}\label{prop_in}
\sup\limits_{0\leq s\leq t\leq\delta}\big\|T(t,s)\big\|_{\cl(\h_\alpha^{m},L^2)}
~=~+\infty,\quad \forall m>0.
\end{equation}
\end{prop}

\section{Nonlinear instability}

In this section, we will prove  the nonlinear ill-posedness result of the three dimensional Prandtl equations stated in Theorem \ref{thm_non}, it will mainly follow the argument of \cite{guo}. First, let us give a preliminary result on the uniqueness of solutions to the linear problem \eqref{lin} as follows.
\begin{lemma}\label{lem_uni}
Suppose that $(u^s,v^s)(t,z)$ is a solution to the problem \eqref{shear} satisfying that
\[\sup\limits_{t\geq0}\Big(\sup\limits_{z\geq0}|(u^s,v^s)|+\int_0^{+\infty}z|(\pd_zu^s,\pd_zv^s)|^2dz\Big)<+\infty.\]
Let $(u,v)\in L^\infty\big(0,T;L^2(\T^2\times\R^+)\big)$ with $(\pd_zu,\pd_zv)\in L^2\big((0,T)\times\T^2\times\R^+\big)$ be a solution to the problem \eqref{lin} with the vanish initial data $(u,v)|_{t=0}=0$. Then, for all $t>0$, $(u,v)\equiv0.$
\end{lemma}
The proof of this lemma is similar to the one given in \cite[Proposition 2.1]{GV-N} or  \cite[Proposition 2.2]{guo}, we omit it here for simplicity.

\vspace{.1in}
\begin{proof}[\bf{Proof of Theorem 2.}]
(1) First, by using \eqref{prop_in}, we know that for the shear flow $(u^s,v^s)$ given in Proposition \ref{prop_uni},  for fixed $\delta_0>0$ and any $m,n\in\N$, there exist $s_n,t_n$ with $0\leq s_n\leq t_n\leq\delta_0$, functions $(u_0^n,v_0^n)(x,y,z)$, and solutions $(u^n_L,v_L^n)$ to the linearized problem \eqref{lin}, such that $(u_L^n,v_L^n)|_{t=s_n}=(u_0^n,v_0^n)$ and
\begin{equation}\label{ass_lin}
\|(u_0^n,v_0^n)\|_{\h_\alpha^m}=1,\quad \|(u_L^n,v_L^n)(t_n)\|_{L^2}\geq n.
\end{equation}

(2) Now, we prove this theorem by contradiction. Assume that the problem \eqref{3dpd} is locally well-posedness for some integer $m\geq0$ in the sense of Definition \ref{def_non}. Denote by
\begin{equation}\label{ass_sh}
(u_{s_n}^s,v_{s_n}^s)(t,z)~:=~(u^s,v^s)(t+s_n,z),
\end{equation}
and
\begin{equation}\label{ass_init}
(u_{0,\delta}^n,v_{0,\delta}^n)(x,y,z)~:=~(u^s,v^s)(s_n,z)+\delta(u_0^n,v_0^n)(x,y,z)
\end{equation}
with a small positive constant $\delta$.
Let $(u_\delta^n,v_\delta^n)$ be the solution to the problem \eqref{3dpd} with the initial data $(u_\delta^n,v_\delta^n)|_{t=0}=(u_{0,\delta}^n,v_{0,\delta}^n)$. Thus, applying Definition \ref{def_non} to two solutions $(u_\delta^n,v_\delta^n)$ and $(u_{s_n}^s,v_{s_n}^s)$, it yields that there exist positive continuous functions $T(\cdot,\cdot)$ and $C(\cdot,\cdot)$, such that
\begin{equation}\label{est_non}\begin{split}
&\|(u_\delta^n-u_{s_n}^s,v_\delta^n-v_{s_n^n}^s)\|_{L^\infty(0,T;L^2(\T^2\times\R^+))}
+\|(u_\delta^n-u_{s_n}^s,v_\delta^n-v_{s_n^n}^s)\|_{L^2(0,T;H^1(\T^2\times\R^+))}\\
&\leq C\delta\|(u_0^n,v_0^n)\|_{\h_\alpha^m}=C\delta,
\end{split}\end{equation}
where $T=T\Big(\|(u_{0,\delta}^n-U_0,v_{0,\delta}^n-V_0)\|_{\h_\alpha^m},\|(u_{s_n}^s-U_0,v_{s_n}^s-V_0)\|_{\h_\alpha^m}\Big)$
and
$C=C\Big(\|(u_{0,\delta}^n-U_0,v_{0,\delta}^n-V_0)\|_{\h_\alpha^m},\|(u_{s_n}^s-U_0,v_{s_n}^s-V_0)\|_{\h_\alpha^m}\Big)$.

Combining \eqref{ass_lin}, \eqref{ass_sh} and \eqref{ass_init}, we know that $\|(u_{0,\delta}^n-U_0,v_{0,\delta}^n-V_0)\|_{\h_\alpha^m}$ and $\|(u_{s_n}^s-U_0,v_{s_n}^s-V_0)\|_{\h_\alpha^m}$ are uniformly bounded in $\delta$ and $n$. Thus, we can take the functions $T(\cdot,\cdot)$ and $C(\cdot,\cdot)$ independent of $\delta$ and $n$, and in the following we use $T,C$ to replace $T(\cdot,\cdot),C(\cdot,\cdot)$ for simplicity.

(3) From the estimate \eqref{est_non}, we know that the sequence $$(\tilde u_\delta^n,\tilde v_\delta^n):=\frac{1}{\delta}(u_\delta^n-u_{s_n}^s,v_\delta^n-v_{s_n^n}^s)$$
is bounded in $L^\infty\big(0,T;L^2(\T^2\times\R^+)\big)\cap L^2\big(0,T;H^1(\T^2\times\R^+)\big)$ uniformly in $\delta$ and $n$, which yields that there is $(u^n,v^n)$ such that, up to a subsequence, as $\delta\rightarrow0,$
\begin{equation}\label{con}
(\tilde u_\delta^n,\tilde v_\delta^n)~\rightharpoonup~(u^n,v^n),\quad{\rm weakly-*~in}\quad L^\infty\big(0,T;L^2(\T^2\times\R^+)\big)\cap L^2\big(0,T;H^1(\T^2\times\R^+)\big),
\end{equation}
and
\begin{equation}\label{est_con}
\|(u^n,v^n)\|_{L^\infty(0,T;L^2(\T^2\times\R^+))}+\|(u^n,v^n)\|_{L^2(0,T;H^1(\T^2\times\R^+))}~\leq~C,\quad\forall n\ge 1.
\end{equation}

Next, since both $(u_\delta^n,v_\delta^n)$ and $(u_{s_n}^s,v_{s_n}^s)$ solve the problem \eqref{3dpd},
we have
\begin{equation}\label{lin_eq}\begin{cases}
\pd_t(\tilde u_\delta^n,\tilde v_\delta^n)+\p_{n}(\tilde u_\delta^n,\tilde v_\delta^n)=\delta N(\tilde u_\delta^n,\tilde v_\delta^n),\\
(\tilde u_\delta^n,\tilde v_\delta^n)|_{t=0}=(u_0^n,v_0^n),
\end{cases}\end{equation}
where $\p_{n}$ is the linearized Prandtl operator at the shear profile $(u_{s_n}^s,v_{s_n}^s)$,
and $N(\cdot, \cdot)$ is the nonlinear term,
\begin{equation}\label{def_n}\begin{split}
N(\tilde u_\delta^n,\tilde v_\delta^n):=\Big(&-\tilde u_\delta^n\pd_x\tilde u_\delta^n-\tilde v_\delta^n\pd_y\tilde u_\delta^n+\int_0^z(\pd_x\tilde u_\delta^n+\pd_y\tilde v_\delta^n)d\tilde z\cdot\pd_z\tilde u_\delta^n,\\
&-\tilde u_\delta^n\pd_x\tilde v_\delta^n-\tilde v_\delta^n\pd_y\tilde v_\delta^n+\int_0^z(\pd_x\tilde u_\delta^n+\pd_y\tilde v_\delta^n)d\tilde z\cdot\pd_z\tilde v_\delta^n\Big).
\end{split}\end{equation}
Therefore, we want to show that the limit function $(u^n,v^n)$ satisfies the linearized Prandtl equations in the sense of distribution. For this, we only need to prove that the right hand side of the equation in \eqref{lin_eq} goes to zero as $\delta\to 0$ in the sense of distribution.

Indeed, note that from \eqref{def_n} the nonlinear term can be rewritten as
\[N(\tilde u_\delta^n,\tilde v_\delta^n)=\Big(N_1(\tilde u_\delta^n,\tilde v_\delta^n),N_2(\tilde u_\delta^n,\tilde v_\delta^n)\Big)\]
with
\[ \begin{split}
&N_1(\tilde u_\delta^n,\tilde v_\delta^n)=-\pd_x(\tilde u_\delta^n)^2-\pd_y(\tilde u_\delta^n\tilde v_\delta^n)+\pd_z\big(\int_0^z(\pd_x\tilde u_\delta^n+\pd_y\tilde v_\delta^n)d\tilde z\cdot\tilde u_\delta^n\big),\\
&N_2(\tilde u_\delta^n,\tilde v_\delta^n)=-\pd_x(\tilde u_\delta^n\tilde v_\delta^n)-\pd_y(\tilde v_\delta^n)^2+\pd_z\big(\int_0^z(\pd_x\tilde u_\delta^n+\pd_y\tilde v_\delta^n)d\tilde z\cdot\tilde v_\delta^n\big).
\end{split}\]
For any compact set $K$ of $[0,T]\times\T^2\times\R^+$ and smooth function $\varphi$ supported in $K$, we have
\[\begin{split}
&\Big|\int_{[0,T]\times\T^2\times\R^+}N_1(\tilde u_\delta^n,\tilde v_\delta^n)\varphi dxdydzdt\Big|\\
&\leq C_{K,\varphi}\int_{K}\Big(|\tilde u_\delta^n|^2+|\tilde u_\delta^n\tilde v_\delta^n|+\big|\int_0^z(\pd_x\tilde u_\delta^n+\pd_y\tilde v_\delta^n)d\tilde z\cdot\tilde u_\delta^n\big|\Big)dxdydzdt\\
&\leq C_{K,\varphi}\Big(\|\tilde u_\delta^n\|_{L^2(K)}^2+\|\tilde v_\delta^n\|_{L^2(K)}^2+\big\|\int_0^z(\pd_x\tilde u_\delta^n+\pd_y\tilde v_\delta^n)d\tilde z\big\|_{L^2(K)}^2\Big),
\end{split}\]
where $C_{K,\varphi}$ is a positive constant depending on $K$ and $W^{1,\infty}$ norm of $\varphi$. From the obvious inequality,
\[\big|\int_0^z(\pd_x\tilde u_\delta^n+\pd_y\tilde v_\delta^n)d\tilde z\big|\leq z^{\frac{1}{2}}\big(\int_{\R_z^+}|\pd_x\tilde u_\delta^n+\pd_y\tilde v_\delta^n|^2dz\big)^{\frac{1}{2}},\]
we get
\[\big\|\int_0^z(\pd_x\tilde u_\delta^n+\pd_y\tilde v_\delta^n)d\tilde z\big\|_{L^2(K)}^2\leq C_K\|\pd_x\tilde u_\delta^n+\pd_y\tilde v_\delta^n\|_{L^2}^2\leq C_K(\|\tilde u_\delta^n\|_{L^2(H^1)}^2+\|\tilde v_\delta^n\|_{L^2(H^1)}^2)\]
for some positive constant $C_K$ depending on $K$. Thus, it follows that
\[
\Big|\int_{[0,T]\times\T^2\times\R^+}N_1(\tilde u_\delta^n,\tilde v_\delta^n)\varphi dxdydz\Big|\leq C_{K,\varphi}(\|\tilde u_\delta^n\|_{L^2(H^1)}^2+\|\tilde v_\delta^n\|_{L^2(H^1)}^2).
\]
Similarly, one can deduce
\[
\Big|\int_{[0,T]\times\T^2\times\R^+}N_2(\tilde u_\delta^n,\tilde v_\delta^n)\varphi dxdydz\Big|\leq C_{K,\varphi}(\|\tilde u_\delta^n\|_{L^2(H^1)}^2+\|\tilde v_\delta^n\|_{L^2(H^1)}^2).
\]
Then, by using the uniform boundedness of $(\tilde u_\delta^n,\tilde v_\delta^n)$ in $L^\infty\big(0,T;L^2(\T^2\times\R^+)\big)\cap L^2\big(0,T;H^1(\T^2\times\R^+)\big)$ with respect to $\delta$, it implies that
the nonlinear term $\delta N(\tilde u_\delta^n,\tilde v_\delta^n)$ converges to zero in the sense of distribution.
Thus, letting $\delta\to 0$ in \eqref{lin_eq}, we obtain that $(u^n,v^n)$ solves the following linear problem in the sense of distribution,
\begin{equation}\label{lin_pr}\begin{cases}
\pd_t(u^n, v^n)+\p_{n}(u^n, v^n)=0,\\
(u^n, v^n)|_{t=0}=(u_0^n,v_0^n).
\end{cases}\end{equation}

(4) Shift the time variable $t$ to $t-s_n$ in \eqref{lin_pr}, and denote by
\[(\tilde u^n,\tilde v^n)(t,\cdot)~:=~(u^n,v^n)(t-s_n,\cdot).\]
Then, \eqref{lin_pr} becomes
\begin{equation}\label{lin_pr1}\begin{cases}
\pd_t(\tilde u^n,\tilde v^n)+\p(\tilde u^n, \tilde v^n)=0,\\
(\tilde u^n,\tilde v^n)|_{t=s_n}=(u_0^n,v_0^n),
\end{cases}\end{equation}
which means that $(\tilde u^n,\tilde v^n)$ solves the linearized problem \eqref{lin} with $(\tilde u^n,\tilde v^n)|_{t=s_n}=(u_0^n,v_0^n)$. By virtue of the uniqueness given in Lemma \ref{lem_uni}, it follows that
$$(\tilde u^n,\tilde v^n)~=~(u_L^n,v_L^n),\quad{\rm on}~[s_n,T].$$
Therefore, from \eqref{ass_lin} and \eqref{est_con} we get a contradiction:
\[
n\leq \|(u_L^n,v_L^n)(t_n)\|_{L^2}=\|(\tilde u^n,\tilde v^n)(t_n)\|_{L^2}=\|(u^n,v^n)(t_n-s_n)\|_{L^2}\leq C,\quad \forall n\ge 1,
\]
as the positive constant $C$, given in \eqref{est_non}, is independent of $n$. So we obtain the proof of Theorem \ref{thm_non}.

\end{proof}

\section{Appendix}

In this Appendix, we present the main steps of the proof of Proposition \ref{prop_2} given in Section 2, which shows that the three dimensional linearized Prandtl equations is well-posed locally in time
when one component of the background tangential velocity, such as $u^s$, is monotonic in the normal variable, and we study the problem in the analytic setting only in the horizontal variable $y$.

\begin{proof}[\bf{Proof of Proposition \ref{prop_2}.}]
Let the solution of the linear problem \eqref{lin} have the form
\[(u,v)(t,x,y,z)~=~\sum\limits_{k\in\cZ}e^{iky}(u_k,v_k)(t,x,z).\]
Plugging this relation into \eqref{lin}, it follows that
\begin{equation}\label{lin_de}\begin{cases}
\pd_tu_k+u^s\pd_xu_k-\pd_z^2u_k+ikv^su_k-u_z^s\int_0^z(\pd_xu_k+ikv_k)d\tilde z=0,\\
\pd_tv_k+u^s\pd_xv_k-\pd_z^2v_k+ikv^sv_k-v_z^s\int_0^z(\pd_xu_k+ikv_k)d\tilde z=0,\\
(u_k,v_k)|_{z=0}=0.
\end{cases}\end{equation}
By assuming that $u^s(t,z)$ is monotonic in $z$, i.e. $\partial_zu^s>0$,  we employ the transformation given in \cite{AWXY} for
the first component of the tangential velocity in the above problem,
\[h_k(t,x,z)~\triangleq~\pd_z\big(\frac{u_k(t,x,z)}{u_z^s(t,z)}\big),~{\rm or}~u_k(t,x,z)~=~u_z^s(t,z)\int_0^zh_k(t,x,\tilde z)d\tilde z,\]
and set
\[\tilde v_k(t,x,z)~\triangleq~\big(v_k-\frac{v_z^s}{u_z^s}u_k\big)(t,x,z).
\]
Then, from the problem \eqref{lin_de} we know that $(h_k,  \tilde v_k)$ satisfies the following problem,
\begin{equation}\label{lin_tr}\begin{cases}
\pd_th_k+u^s\pd_xh_k-\pd_z^2h_k-2\pd_z\big(\frac{u_{zz}^s}{u_z^s}h_k\big)+ik(v^sh_k-\tilde v_k)=0,\\
\pd_t\tilde v_k+u^s\pd_x\tilde v_k-\pd_z^2\tilde v_k+ikv^s\tilde v_k-2u_z^s\pd_z(\frac{v_z^s}{u_z^s})h_k=0,\\
\pd_zh_k|_{z=0}=0,\quad \tilde v_k|_{z=0}=0,
\end{cases}\end{equation}
where we use $\pd_z^2u_k|_{z=0},~u_{zz}^s|_{z=0}=0$ to derive the boundary condition of $h_k$.

For the problem \eqref{lin_tr}, by the energy method
one can have
\begin{equation}\label{est_tr}\begin{split}
&\quad\|(h_k,\tilde v_k)(t,\cdot)\|_{K_\alpha^{m}}^2+\int_0^t\|(\pd_zh_k,\pd_z\tilde v_k)(s,\cdot)\|_{K_\alpha^{m}}^2ds\\
&\leq C\Big(\|(h_k,\tilde v_k)(0,\cdot)\|_{K_\alpha^{m}}^2+\max\{1,|k|\}\int_0^t\|(h_k,\tilde v_k)(s,\cdot)\|_{K_\alpha^{m}}^2ds\Big),
\end{split}\end{equation}
where the positive constant $C$ depends on $\alpha$ and $(u^s,v^s)$. Applying the Gronwall inequality to \eqref{est_tr}, it implies that there exists a $\rho>0$, depending on $\alpha$ and $(u^s,v^s)$, such that
\begin{equation}\label{est_tr1}\|(h_k,\tilde v_k)(t,\cdot)\|_{K_\alpha^{m}}^2+\int_0^t\|(\pd_zh_k,\pd_z\tilde v_k)(s,\cdot)\|_{K_\alpha^{m,0}}^2ds\leq Ce^{\rho |k|t}\|(h_k,\tilde v_k)(0,\cdot)\|_{K_\alpha^{m}}^2.
\end{equation}

From the assumption \eqref{ass_ini}, we have
\begin{equation}\label{ini_tr}
\|(h_k,\tilde v_k)(0,\cdot)\|_{K_\alpha^m}~\leq~C_0e^{-\beta |k|},
\end{equation}
for some positive constant $C_0$. As $u^s\in C(\R^+;W_\alpha^{3,\infty}(\R^+))$, one has
\begin{equation}\label{est_tr2}
\|u_k(t,\cdot)\|_{K_\alpha^m}=\|(u_z^s\int_0^zh_kd\tilde z)(t,\cdot)\|_{K_\alpha^m}\leq C_1\|h_k(t,\cdot)\|_{K_\alpha^m},
\end{equation}
with the constant $C_1>0$ depending on $\alpha$ and $u^s$. Therefore, from the estimates  \eqref{est_tr1}-\eqref{est_tr2} and the relation $v_k=\tilde v_k+\frac{v_z^s}{u_z^s}u_k$ it follows that
\begin{equation}\label{est_de}
\|(u_k, v_k)(t,\cdot)\|_{K_\alpha^m}~\leq~C_2e^{-(\beta-\rho t)|k|},
\end{equation}
where the positive constant $C_2$ is independent of $k$. From the estimate \eqref{est_de} we complete the proof of
this proposition.
\end{proof}

\vspace{.15in}

{\bf Acknowledgements :}
The first two authors' research was supported in part by
National Natural Science Foundation of China (NNSFC) under Grants
No. 10971134, No. 11031001 and No. 91230102. The last author's research was supported by the General Research Fund of Hong Kong,
CityU No. 103713.

\end{document}